\documentclass[12pt]{amsart}
\usepackage[T1]{fontenc}
\usepackage[utf8]{inputenc}
\usepackage{color,amssymb,amsmath,mathrsfs,enumerate,esint}
\usepackage{a4wide}
\usepackage{hyperref}
\usepackage{pgf,tikz,pgfplots}
\pgfplotsset{compat=1.14}
\usepackage{mathrsfs}
\usetikzlibrary{arrows}
\theoremstyle{plain}
\newtheorem{Theorem}{Theorem}[section]
\newtheorem{Lemma}[Theorem]{Lemma}
\newtheorem{Question}[Theorem]{Question}

\newtheorem{Corollary}[Theorem]{Corollary}
\newtheorem{Example}[Theorem]{Example}
\theoremstyle{definition}
\newtheorem{Remark}[Theorem]{Remark}

\DeclareMathOperator{\er}{\mathbb{R}}
\DeclareMathOperator*{\dee}{d}
\DeclareMathOperator*{\loc}{loc}
\DeclareMathOperator{\supp}{supp}
\title[Optimal G-N inequality for function spaces]{Optimal Gagliardo--Nirenberg interpolation inequality for rearrangement invariant spaces}
\author{Karol Le\'snik}

\author{Tom\' a\v s Roskovec}

\author{Filip Soudský}

\address{K.~Le\'snik: Faculty of Mathematics and Computer Science,  Adam Mickiewicz University in Pozna{\'n},
  Uniwersytetu Pozna{\'n}skiego 4,
  61-614 Pozna{\'n}, Poland}
\email{klesnik@vp.pl}

\address{T.~G.~Roskovec: Faculty of Economics, University of South Bohemia, Studentsk\' a 13, \v Cesk\' e Bud\v ejovice, Czech Republic; \newline
Faculty of Science, University of Hradec Kr\' alov\' e, Hradeck\' a 1285, Hradec Kr\' alov\' e }
\email{troskovec@ef.jcu.cz}

\address{F.~Soudsk\' y: Faculty of Science, Humanities and Education, Technical University of Liberec, Studentská 1402/2, Liberec, Czech Republic; \newline
Faculty of Science, University of Hradec Kr\' alov\' e, Hradeck\' a 1285, Hradec Kr\' alov\' e }
\email{filip.soudsky@tul.cz}

\begin{document}

\begin{abstract}

    We prove optimality of the Gagliardo-Nirenberg inequality
    \[
    \|\nabla u\|_{X}\lesssim\|\nabla^2 u\|_Y^{1/2}\|u\|_Z^{1/2},
    \]
    where $Y, Z$ are rearrangement invariant Banach function spaces and $X=Y^{1/2}Z^{1/2}$ is the Calder\'on--Lozanovskii space. By optimality we mean that for a certain pair of spaces on the right-hand side, one cannot reduce the space on the left-hand, remaining in the class of rearrangement invariant spaces. The optimality for the Lorentz and Orlicz spaces is given as a consequence, exceeding previous results. We also discuss pointwise inequalities, their importance and counterexample prohibiting an improvement.
\end{abstract}
\maketitle
\begin{flushright}
\it Dedicated to the memory of \\
Professor Jan Malý \\
(1955--2021)
\end{flushright}

\section{Introduction}
The Sobolev-Gagliardo--Nirenberg interpolation inequality is typically inequality in the form
\begin{equation}\label{eq:GNgeneral}
\|\nabla^j u\|_{X}\lesssim\|\nabla^k u\|_Y^{\theta}\|u\|_Z^{1-\theta},
\end{equation}
where $0\leq j<k$ are the degrees of derivatives, $\theta\in[j/k, 1]$ and $X, Y, Z$ are function spaces. The general form is the product of continued development. The first attempts on estimating an intermediate derivative were made at the beginning of the 20th century by Landau (see \cite{Lan}) who proved inequality
$$
\| u'\|_{L^\infty}\lesssim\| u''\|_{L^\infty}^{1/2}\|u\|_{L^\infty}^{1/2}.
$$
Later on, Kolmogoroff proved the result for higher-order derivatives \cite{Kol}. In 1958 the problem was revisited by Nash (see \cite{N}) and Ladyzenskaya who proved
$$
\|u\|_{L^{4}}\lesssim\|\nabla u\|_{L^2}^{1/2}\|u\|_{L^2}^{1/2}
$$
in dimension two and also similar type result in dimension three. Her research \cite{Lady} was motivated by the study of mathematical models for hydro-mechanics. Later, Stein in \cite{STEIN} improved Landau-Kolmogoroff's result in one dimension as
$$
\|u'\|_{L^p}\lesssim\|u''\|^{1/2}_{L^p}\|u\|^{1/2}_{L^p}
$$
for $p\in[1,\infty)$.
The real brake-through came in 1959 as Gagliardo \cite{Ga} and Nirenberg \cite{Ni} independently introduced the general version of Gagliardo--Nirenberg interpolation inequality for Lebesgue spaces which took the form
\begin{equation}\label{GNCI}
\|\nabla^j u\|_{L^p}\lesssim\|\nabla^k u\|^\theta_{L^r}\|u\|_{L^q}^{1-\theta}.
\end{equation}
The formula holds for parameters $0\leq j<k$, $\theta=j/k$, and $p,q,r$ satisfying $k/p=j/r+(k-j)/q$. Also a slightly more complicated setting for the case $\theta\in[j/k,1]$ is proven therein. The modern versions of the proof can be found in \cite{FIFOROSO2} or \cite{Leon}.

Nowadays, the field is too wide for being covered without omitting a lot. The primary motivation is the boundedness of the norm guaranteeing desired property for the solution of PDE, see, for example, \cite{Bocc, DelPi, PePo}. As the methods and the function spaces scales refine, the demand for a finer version of such an inequality appears. As for other applications, we may mention the chain rule for Sobolev spaces \cite{MoRoSo} or the study of the boundedness of bilinear multipliers \cite{Grafakos, Slav}. The problem is developed in several branches, and the term Gagliardo--Nirenberg inequality may be used for several slightly different settings.

We may split the field concerning \eqref{eq:GNgeneral} based on value of $j$. Either $j>0$ as by Gagliardo or $j=0$ as by Ladyzhenskaya.  Interesting modern results for $j>0$ are for example setting the inequality in the Orlicz spaces scale \cite{KP} or relaxing $L^{\infty}$ space into BMO in \cite{SP}. In literature the more dominant version is corresponding to Ladyzhenskaya's result, where \eqref{eq:GNgeneral} is set for $j=0$ and typically $k=1$. This version is studied in \cite{MRR, DelPi} and used in above mentioned \cite{Bocc, PePo}.

The generalizations have been done in several directions. Even in the original paper, Nirenberg extends the Lebesgue scale into negative $p$ and claims the result for H\"older spaces. The other direction is to extend $j$ and $k$ from integer number into real, more precisely using the fractional Sobolev space. Such a setting may be found in \cite{BM, BreMir}.

Yet another possible modification is the so-called nonlinear or strongly nonlinear Gagliardo--Nirenberg inequality
$$
\|\nabla^j u\|_{X}\lesssim F(\nabla^k u,u),
$$
where $F$ is a functional, such as in \cite{RiSt, KaPe}. These types may be beneficial in special cases of applications where standard multiplicative versions may not be enough.

To map the development, it is worth mentioning the involvement of modern tools in both formulation and the proving methods. We already mention fractional Sobolev spaces, which not only create new settings but also help to reprove older results with ease \cite{BM}. The problem was also investigated in the Lorentz spaces scale. Relaxations from the Lebesgue spaces to the Lorentz spaces scale were given for $j=0$ in \cite{MRR} and, recently, in \cite{wiki}. See also \cite{KCh} for variable Lebesgue spaces setting and \cite{BM, BSP1, BSP2} for fractional Sobolev and Besov spaces setting.


\section{Results}
 We start the paper by reformulating (or reproving) the version of Gagliardo--Nirenberg inequality for rearrangement invariant Banach function spaces stated already in \cite{FiFoRoSo} (cf. \cite{KCh}) in another form. We do so, since the statement from \cite{FiFoRoSo}  allows only very limited analysis of optimality, as we explain in Section \ref{subsec.opakovanimatkamoudrosti}. The rougher versions of optimality may be found in \cite{FiFoRoSo} or even \cite{Ni}, where the optimality from the point of view of fundamental functions is given, covering for example the optimality on the scale of Lebesgue and Orlicz spaces. Some other special cases may be found elsewhere in the literature.

\begin{Theorem}[Gagliardo--Nirenberg inequality for r.i.B.f. spaces]\label{GBFSR}
Let $Y,Z$ be rearrangement invariant Banach function spaces and let $X=Y^{j/k}Z^{1-j/k}$ be defined by the Calder\'on-Lozanovskii construction, where $1\leq j< k$ and upper Boyd indices of both spaces $Y$ and $Z$ are smaller than 1. Then
\begin{equation}\label{ME}
\|\nabla^j u\|_{X}\lesssim\|\nabla^k u\|_Y^{\frac{j}{k}}\|u\|_Z^{1-\frac{j}{k}}
\end{equation}
holds for all $u\in W^{k,1}_{\loc}(\mathbb{R}^d)\cap Z$. 
\end{Theorem}
The essential contribution is the usage of the maximal operator and producing the pointwise version of the Gagliardo--Nirenberg inequality. The first article in this direction is probably \cite{KufMa} by Kufner and Maz$'$ya, but the stronger very useful result was first proven by Ka{\l}amajska \cite{Kal}. 
We use the following weaker version of Ka{\l}amajska result 
\begin{equation}\label{eq:pointwise}
|\nabla^j u|\lesssim (M\nabla^k u)^{\frac{j}{k}}(Mu)^{1-\frac{j}{k}},\ \ {\rm\ where\ }0<j\leq k
\end{equation}
proved later independently by Maz$'$ya and Shaposchnikova in \cite{MS}.

It is worth mentioning here, that pointwise estimates have already been applied to prove more general forms of the Gagliardo--Nirenberg inequality, such as the Gagliardo--Nirenberg for Orlicz spaces considered by Ka{\l}amajska, Pietruska-Pa{\l}uba, and Krbec (see \cite{KP} and \cite{KM}) or even for general Banach function spaces \cite{KCh, FiFoRoSo}.

Later, Lokharu focused on different types of maximal operators in \cite{Lok1, Lok2}.  Notice also that early version of Theorem \ref{GBFSR} appeared in the literature before \cite{FiFoRoSo} in \cite{KCh}, but without the focus on the optimality.

Having Theorem \ref{GBFSR} in suitable form, we are ready to discuss optimality of it in the setting of rearrangement invariant Banach function spaces. By optimality we mean, that having fixed two spaces, say $Y, Z$, there is no space $B$ such that it is smaller than $X$ and \eqref{ME} still holds.

It is already known that Theorem \ref{GBFSR} is sharp in the scale of Lebesgue space. It follows from the so-called scaling argument, see \cite[Thm. 1.1]{FiFoRoSo}. 
In fact, the scaling argument applied to the setting of rearrangement invariant Banach function spaces gives a  necessary condition on the fundamental function of space $X$ satisfying \eqref{ME}, as was proved in \cite{FiFoRoSo}.  \eqref{eq:GNgeneral}. 
However, (except the limit cases) there is a whole universe of rearrangement invariant Banach function spaces with the same fundamental function. Thus to decide whether the choice of $X=Y^{j/k}Z^{1-j/k}$ in \eqref{ME} is optimal in the scale of rearrangement invariant Banach function spaces we need to introduce an essentially different and more sensitive approach, than the scaling argument. We fix our attention only on the most classical case of $j=1,k=2$. This approach is the merit of the following optimality theorem.

\begin{Theorem}[Optimality of Gagliardo--Nirenberg inequality for r.i.B.f. spaces]\label{optimal}
Let $Y,Z$ be r.i.B.f. spaces and $X=Y^{1/2}Z^{1/2}$. Assume that the Gagliardo--Nirenberg inequality
\begin{equation*}
\|\nabla u\|_{B}\lesssim\|\nabla^2 u\|_Y^{\frac{1}{2}}\|u\|_Z^{\frac{1}{2}} 
\end{equation*}
holds for all $u\in W^{2,1}_{\loc}(\mathbb{R}^d)\cap Z$ and some rearrangement invariant Banach function space $B$.  
\begin{enumerate}[\upshape(i)]
\item If $Z\cap L^{\infty}\subset Y\cap L^{\infty}$, then $X\cap L^{\infty}\subset B\cap L^{\infty}$.
\item
If $Z\cap L^{1} \subset Y\cap L^{1}$, then $X\cap L^{1} \subset B\cap L^{1}$.
\end{enumerate}
\end{Theorem}

\begin{Remark}
Notice that we do not need restrictions on upper Boyd indices of $Y$ and $Z$, in contrast to Theorem \ref{GBFSR}. 
\end{Remark} 

When both assumptions of points (i) and (ii) are satisfied at the same time, we get the complete optimality of $X=Y^{j/k}Z^{1-j/k}$. 

\begin{Corollary}
Let $Y,Z$ be r.i.B.f. spaces satisfying $Z\subset Y$. If 
the Gagliardo--Nirenberg inequality
\begin{equation*}
\|\nabla u\|_{B}\lesssim\|\nabla^2 u\|_Y^{\frac{1}{2}}\|u\|_Z^{\frac{1}{2}}
\end{equation*}
holds for all $u\in W^{2,1}_{\loc}(\mathbb{R}^d)\cap Z$ and some r.i.B.f. space $B$, then $Y^{1/2}Z^{1/2}\subset B$.
\end{Corollary}

Concluding, we see that the choice of $X=Y^{1/2}Z^{1/2}$ is optimal among all r.i.B.f. spaces in the Gagliardo--Nirenberg inequality  \eqref{ME}, provided that $Z\subset Y$. However, the assumption $Z\subset Y$ is quite restrictive and does not apply to the most classical r.i. spaces (Lebesgue spaces, Lorentz spaces, etc.), since usually there is no inclusion between such spaces over $\mathbb{R}_+$. It appears, however, that manoeuvring between points (i) and  (ii) of Theorem \ref{optimal} we can use it to give an almost complete answer to the question about the optimality of  \eqref{ME} among Lorentz spaces posted in \cite[Remark 2.7]{FiFoRoSo}.

\begin{Corollary}\label{Lorentz}
Let $1< R,Q< \infty$ and $1\leq r,q\leq \infty$. Then  the  Gagliardo--Nirenberg inequality
\begin{equation}\label{eq>GNlor0}
\|\nabla u\|_{P,p}\lesssim\|\nabla^2 u\|_{Q,q}^{\frac{1}{2}}\|u\|_{R,r}^{\frac{1}{2}}
\end{equation}
holds for all $u\in W^{2,1}_{\loc}(\mathbb{R}^d)\cap L^{R,r}$ with $P=(1/(2Q)+1/(2R))^{-1}$ and $p= (1/(2q)+1/(2r))^{-1}$.

On the other hand, if one of the following conditions holds 
\begin{enumerate}[(i)]
    \item $R\not =Q$, 
    \item $R=Q$ and $r<q$,
\end{enumerate}
then the  Gagliardo--Nirenberg inequality
\begin{equation}\label{eq>GNlor}
\|\nabla u\|_{P',p'}\lesssim\|\nabla^2 u\|_{Q,q}^{\frac{1}{2}}\|u\|_{R,r}^{\frac{1}{2}}
\end{equation}
holds for all $u\in W^{2,1}_{\loc}(\mathbb{R}^d)\cap L^{R,r}$ if and only if  $P'=P$ and $p'\geq p$.
\end{Corollary}

Notice, that the corollary above does not cover the case of $Q=R$ and $r>q$. 
\begin{Question}
Is the Gagliardo--Nirenberg inequality  \eqref{eq>GNlor} valid for some $p'<p:=(1/(2q)+1/(2r))^{-1}$ in case $1< P'=R=Q<\infty$ and $1\leq q<r$?
\end{Question}

Our main results give also new information about the Gagliardo--Nirenberg inequality for Orlicz spaces. 

\begin{Corollary}\label{Orlicz}
Let $\varphi,\varphi_1,\varphi_2$ be Young functions such that $\varphi^{-1}\approx \sqrt{\varphi_1^{-1}\varphi_2^{-1}}$ and upper Boyd indices of $L^{\varphi_1}$ and $L^{\varphi_2}$ are smaller than 1.  Then  the  Gagliardo--Nirenberg inequality 
\begin{equation}\label{GNOS0}
\|\nabla u\|_{\varphi}\lesssim \|\nabla^2 u\|_{\varphi_1}^{\frac{1}{2}}\|u\|_{\varphi_2}^{\frac{1}{2}}
\end{equation}
holds for all $u\in W^{2,1}_{\loc}(\mathbb{R}^d)\cap L^{\varphi_2}$. 
On the other hand, if $L^{\varphi_2}\subset L^{\varphi_1}$ and the Gagliardo--Nirenberg inequality 
\begin{equation}\label{GNOS1}
\|\nabla u\|_{B}\lesssim \|\nabla^2 u\|_{\varphi_1}^{\frac{1}{2}}\|u\|_{\varphi_2}^{\frac{1}{2}}
\end{equation}
holds for all $u\in W^{2,1}_{\loc}(\mathbb{R}^d)\cap L^{\varphi_2}$ and some r.i.B.f. space $B$, then $L^{\varphi}\subset B$. 
\end{Corollary}

Notice that \eqref{GNOS0} has been already proved in \cite{FiFoRoSo} (cf. \cite{KP}). On the other hand, it was shown therein that $L^{\varphi}$ as above is optimal among all Orlicz spaces, while our result says that it is optimal also in the class of r.i.B.f. spaces, provided $L^{\varphi_2}\subset L^{\varphi_1}$.

\section{Preliminaries}
Given some positive-valued  $F(u), G(u)$ we write
$$
F(u)\lesssim G(u)
$$
if there exists a constant $C>0$ (independent on $u$) such that
$$
F(u)\leq C G(u).
$$
If both
$F(u)\lesssim G(u)$ and $G(u)\lesssim F(u)$, then we write
$$
F(u)\approx G(u).
$$

\subsection{Function spaces instruments} 
In our paper, we shall use the standard notation. Given a Lebesgue-measurable set $A\subset\mathbb{R}^d$, we shall denote its Lebesgue measure by $|A|$.
Symbol $L^0(A)$ denote the space of all measurable functions over the given set $A$ that are finite almost everywhere. 

Given a Banach space $X\subset L^0$ we will call it a {\it Banach function space} (B.f. space, for short) if it has the following properties:
\begin{itemize}
\item ideal property, i.e. given $f\in L^0$ and $g\in X$ with $|f|\leq g$ there holds $\|f\|_X\leq \|g\|_X$. 
\item Fatou property, i.e. given $f\in L^0$ and an increasing sequence $(f_n)\in X$ such that $0\leq f_n\to f$ a.e. and $\sup_n\|f_n\|_X<\infty$, there holds $f\in X$ and $\|f\|_X=\sup\limits_n\|f_n\|_X$. 
\end{itemize}

Given a measurable function $f\in L^0(\mathbb{R}^d)$ and $\alpha\in\mathbb{R}$ we shall denote its \textit{level set} shortly by
$$
\{f>\alpha\}:=\{x\in\mathbb{R}^d: f(x)>\alpha\}
$$
and similarly in case of $\{f\geq\alpha\}, \{f<\alpha\}, \{f\leq\alpha\}$.
Throughout the paper, a \textit{distribution function of f} is denoted by
$$
f_*(t):=\left|\{|f|>t\}\right|
$$
and the \textit{non-increasing rearrangement} by
$$
f^*(t):=\inf_{s}\{f_*(s)\leq t\}.
$$

A B.f. space $X$ is called \textit{rearrangement invariant B.f. space} (we use the abbreviation r.i.B.f. space) if given two functions $f\in L^0$ and $g\in X$ which are equidistributed (i.e. $f_*=g_*$), there holds $f\in X$ and $\|f\|_X=\|g\|_X$. 

To simplify the notion, we make use of the \textit{Luxemburg representation theorem}. It says that for each r.i.B.f. space $X$ over arbitrary $A\subset \mathbb{R}^d$ there is an r.i.B.f. space  $\overline{X}$ over $(0,\infty)$ with the Lebesgue measure, 
such that for every measurable $f:A\to \mathbb{R}$ one has
\begin{equation}\label{eq:Lux-rep}
\|f\|_X=\|f^*\|_{\overline{X}}.
\end{equation}
Thanks to it, we can and we do restrict our attention to r.i.B.f. spaces defined over $(0,\infty)$ and do not consider r.i.B.f. spaces over other domains. This is, henceforth, speaking about r.i.B.f. space $X$ we always mean that it is defined over $(0,\infty)$. Yet, we work with functions defined over $\mathbb{R}^d$. Consequently, having such a function and an r.i.B.f. space $X$, by writing 
$\|f\|_{X}$
we mean just 
$ \|f^*\|_{X}$.

Given a r.i.B.f. space $X$ we define
$$
X_{\loc}:=\{u\in L^0: u^{*}\chi_{(0,1)}\in \overline{X}\}.
$$
Moreover, we write
$$
X\stackrel{\loc}{\hookrightarrow}Y,
$$
if there exists a constant $C>0$ such that
$$
\|f^*\chi_{(0,1)}\|_{\overline{X}}\leq C\|f^*\chi_{(0,1)}\|_{\overline{Y}},\quad \forall f\in L^{0}.
$$
We also need the \textit{Hardy operator}, which is defined for $u\in L^{1}_{\loc}$ by the formula
$$
Au(t):=\frac{1}{t}\int_0^t u(s)\textup{d}s, \ \ t>0.
$$
Further, the symbol $u^{**}$ will be used for the \textit{Hardy-Littlewood maximal function} of $u$ defined by
$$
u^{**}(t):=A(u^*)(t), \ \ t>0.
$$
It is known, that each r.i.B.f. space $X$ satisfies the \textit{Hardy-Littlewood-Polya principle} (or the \textit{majorant property}), which states that for $g\in X$ and $f\in L^{1}_{\loc}$
\begin{equation}\label{eq:HLP}
\textup{if }f^{**}(t)\leq g^{**}(t) \text{ for all } t>0 \textup{, then\ } f\in X {\rm \ and\ }\|f\|_X\leq\|g\|_X.
\end{equation}

The maximal function defined above is strictly connected with the \textit{maximal operator} which is defined for any function locally integrable on $\mathbb{R}^d$, by the formula
$$
Mu(x):=\sup_{Q\ni x}\frac{1}{|Q|}\int_Q |u(y)|\textup{d} y,\ \ x\in \mathbb{R}^d,
$$
where the supremum on the right-hand side is taken over all $d$--dimensional cubes $Q$ containing $x$. 
The \textit{Riesz-Herz equivalence} states that
\begin{equation}\label{eq:Riesz-Herz}
(Mf)^*(t)\approx f^{**}(t)\quad \text{ for all } t\in(0,\infty),
\end{equation}
where the constant of equivalence depends only on the dimension $d$ (see \cite[Theorem 3.8]{BS}). The \textit{dilation operator} $D_s$ is defined for $s>0$ by the formula
$$
D_sf(t):=f(st).
$$
For more information on r.i.B.f. spaces we refer for instance to books \cite{BS, FS,KPS}.

The most significant examples of r.i.B.f. spaces except for the Lebesgue spaces are Lorentz and Orlicz spaces. Let us recall its definitions. 
For $1\leq P,p\leq\infty$ the {\rm Lorentz space} $L^{P,p}$ is given by the (quasi--) norm
$$
\|f\|_{P,p}=\left(\int_0^{\infty}[t^{1/P}f^*(t)]^p\frac{dt}{t}\right)^{1/p}.
$$
In the case $1<P$ and $1\leq p\leq \infty$ the functional $\|\cdot \|_{P,p}$ is already equivalent to the norm, thus we will treat $L^{P,p}$ as a r.i.B.f. space. When $P=1$ we consider only $p=1$, since otherwise $L^{1,p}$  is no more a B.f. space.

A continuous, convex, increasing function $\varphi:[0,\infty)\to [0,\infty)$ satisfying $\varphi(0)=0$ is called the {\it Young function}. For a given Young function $\varphi$ we define the Orlicz space $L^{\varphi}$ by
\[
L^{\varphi}=\{f\in L^0:\int_0^{\infty}\varphi\left(\frac{f(t)}{\lambda}\right)dt <\infty {\rm \ for\ some \ }\lambda>0 \}
\]
with the Luxemburg norm
\[
\|f\|_{\varphi}=\inf \{\lambda>0 :\int_0^{\infty}\varphi\left(\frac{f(t)}{\lambda}\right)dt \leq 1\}.
\]

\subsection{Derivatives}
To denote the distributional derivative in one-dimensional case we use symbols $u',u'',u^{(k)}$. The symbols $\nabla,\nabla^k$ denote the distributional gradient, respectively the distributional gradient of order $k$. We set
$$
\left|\nabla^k u\right|:=\sum_{|\alpha|=k}\left|\frac{\partial^\alpha u}{\partial x^\alpha}\right|.
$$
Given a r.i.B.f. space $X$ we set 
$$
\|\nabla^k u\|_X:=\||\nabla^k u|\|_X.
$$
We also define the maximal operator of the $k$-order gradient by
$$
M(\nabla^k u)(x):=\sup_{Q\ni x}\frac{1}{|Q|}\int_{Q}|\nabla^k u|\textup{d}x.
$$

Notion $W^{k,p}(A)$ will be used for the space of functions whose distributional derivatives up to degree $k$ belong to $L^{p}(A)$, sub-index ``loc'' is added if the distributional derivatives belong only to $L^{p}_{\loc}(A)$.

\subsection{Calder\' on--Lozanovskii construction and pointwise product spaces}

Given $0<\theta<1$ and  two  B.f. spaces $X, Y$ over the same measure space, the {\it Calder\'on--Lozanovskii space} $X^{\theta}Y^{1-\theta}$ is defined as 
\[
X^{\theta}Y^{1-\theta}=\{f\in L^0 :|f|\leq g^{\theta}h^{1-\theta} {\rm \ for \ some\ }0\leq g\in X,0\leq h\in Y \}
\]
and equipped with the norm 
\[
\|f\|_{X^{\theta}Y^{1-\theta}}=\inf\left\{ \max\left\{\|g\|_X,\|h\|_Y\right\}:|f|\leq g^{\theta}h^{1-\theta}, 0\leq g\in X,0\leq h\in Y\right\}.
\]
It is easily seen that one can replace ``$\leq$'' by ``$=$'' in the above definitions. Moreover, the following inequality holds true for each $ g\in X,h\in Y$ (see \cite{Mal89})
\begin{equation}\label{complex}
\| |g|^{\theta}|h|^{1-\theta}\|_{X^{\theta}Y^{1-\theta}}\leq \|g\|_X^{\theta}\|h\|_Y^{1-\theta}.
\end{equation} 

Given a B.f. space $X$ and $\alpha>1$, we define the \textit{$\alpha$-convexification of  $X$} (respectively, \textit{$\alpha$-concavification} when $0<\alpha<1$) by
$$
X^\alpha:=\{f\in L^0:|f|^\alpha\in X\}
$$
with the (quasi--) norm given by
$$
\|f\|_{X^{\alpha}}:=\left(\||f|^\alpha\|_X\right)^{1/\alpha}.
$$
It is easy to see that convexification is just a special case of the Calder\'on--Lozanovskii construction, namely
$$
X^\alpha=X^{\frac{1}{\alpha}}(L^{\infty})^{1-\frac{1}{\alpha}}.
$$

In Section \ref{Optimality}, we will need once more constructions intimately connected with Calder\'on--Lozanovskii spaces. Given two B.f. spaces $X, Y$ over the same measure space we define their pointwise product $X\odot Y$ by 
\[
X\odot Y=\{f\in L^0:|f|= gh  {\rm \ for \ some\ } g\in X,h\in Y\},
\]
equipped with the quasi-norm 
\[
\|f\|_{X\odot Y}=\inf\{ \|g\|_X\|h\|_Y:|f|= gh  {\rm \ for \ some\ }g\in X,h\in Y\}.
\]
We will need the following identification from \cite{KLM14}
\begin{equation}\label{propCL1}
X^{\theta}Y^{1-\theta}=X^{\frac{1}{\theta}}\odot Y^{\frac{1}{1-\theta}}.
\end{equation}


\section{New shape of Gagliardo--Nirenberg inequality for r.i.B.f. spaces}\label{subsec.opakovanimatkamoudrosti}
We start this section with the proof of Theorem \ref{GBFSR}, which actually is an immediate consequence of \eqref{eq:pointwise}.

\proof[Proof of Theorem \ref{GBFSR}]
Let $u\in W^{k,1}_{\loc}(\mathbb{R}^d)$ and $1\leq j<k$. We have by \eqref{eq:pointwise}
$$
\begin{aligned}
\|\nabla^j u\|_X&\lesssim \|(M(\nabla^k u))^{\frac{j}{k}} (Mu)^{1-\frac{j}{k}}\|_X.
\end{aligned}
$$
Then by properties of the Calderon--Lozanovski construction
$$
\begin{aligned}
\|(M(\nabla^k u))^{\frac{j}{k}} (Mu)^{1-\frac{j}{k}}\|_X
\leq \|M(\nabla^{k}u)\|_{Y}^{\frac{j}{k}}\|Mu\|_{Z}^{\frac{k-j}{k}}.
\end{aligned}
$$
Finally, we replace the maximal operators of functions by functions themselves in both norms on the right-hand side. The boundedness of the maximal operator in considered spaces (based on the upper Boyd index being smaller than 1) gives
$$
\|\nabla^j u\|_X\lesssim  \|(\nabla^{k}u)\|_{Y}^{\frac{j}{k}}\|u\|_{Z}^{\frac{k-j}{k}},
$$
as desired.
\endproof

Let us note that without the assumption about the boundedness of the maximal operator, we cannot obtain the last estimate.
Neither counterexamples nor the proof of the validity of theorem without the assumption on the Boyd indices is known. Let us however notice, that the original proof of the Gagliardo--Nirenberg inequality for Lebesgue spaces allows also $L^1$ spaces (cf. \cite{Ga, Ni, FIFOROSO2}),  the space where the maximal operator is not bounded.
On the other hand, proofs of the inequality for Orlicz, Lorentz or more general spaces (up to our knowledge) are based on the same kind of pointwise inequalities, while such an approach requires assumptions on Boyd indices (or generally on boundedness of maximal operator). 
\begin{Question}\label{Q1}
Does Theorem \ref{GBFSR} hold without the assumption on Boyd indices of $Y$ and $Z$? 
\end{Question}

One of the possible ways of removing assumptions on Boyd indices would be a slight improvement of \eqref{eq:pointwise}.
In fact, rewriting \eqref{eq:pointwise}  in the form
$$
(\nabla^j u)^{**}\lesssim (|\nabla^k u|^{**})^{\frac{j}{k}}(|u|^{**})^{1-\frac{j}{k}},
$$
by the Riesz-Herz equivalence, we observe that it would be sufficient for our purpose if we could put the geometric mean of $|u|$ and $|\nabla^2 u|$  inside the double-star operator. Precisely, the question is if the following inequality holds
\begin{equation}\label{false}
(\nabla^j u)^{**}\lesssim (|\nabla^k u|^{\frac{j}{k}}|u|^{1-\frac{j}{k}})^{**}?
\end{equation}
This inequality would then, by the simple use of Hardy-Littlewood-Polya principle, imply that Theorem \ref{GBFSR} holds with no restriction on Boyd indices. However, the counterexample below shows that one cannot hope for \eqref{false}.

\begin{Example}\label{Example Mount Filip}
Let $u_n:\er\to \er$ be a sequence of functions defined in the following way
$$
u_n(0):=0, \quad u_n'(0)=0\quad u_n''(s):=n\left(\chi_{(0,\frac{1}{n})}-\chi_{(1-\frac{1}{n},1+\frac{1}{n})}+\chi_{(2-\frac{1}{n},2)}\right).
$$
One easily verifies that all the functions are supported in interval $[0,2]$ and $u_n(2)=u_n'(2)=0$. Moreover, we estimate
$$
(u_n')^*(s)=\chi_{(0,2-\frac{4}{n})}+\left(1-n\left(s-\left(2-\frac{4}{n}\right)\right)\right)\chi_{(2-\frac{4}{n},2)},
$$
hence we have
$$
1\geq(u_n')^{**}(2)\geq\frac{1}{2}\int_0^2 \chi_{(0,2-\frac{4}{n})}=\frac{2-\frac{4}{n}}{2}=1-\frac{2}{n}
$$
On the other hand
$$
\begin{aligned}
(|u_n''|^{\frac{1}{2}}|u_n|^{\frac{1}{2}})^{**}(2)&=\frac{1}{2}\int_0^2|u_n''|^{\frac{1}{2}}|u_n|^{\frac{1}{2}}(s)\dee s\\
&\leq\frac{1}{2}\left(\int_0^{\frac{1}{n}}\sqrt{n}\dee s+\int_{1-\frac{1}{n}}^{1+\frac{1}{n}}\sqrt{n}\dee s+\int_{2-\frac{1}{n}}^1\sqrt{n}\right)\\
&=\frac{2}{\sqrt{n}}
\end{aligned}
$$
One readily checks that
$$
\lim_{n\to\infty} (u_n')^{**}(2)=1
$$
while
$$
\lim_{n\to\infty} (|u_n''|^{\frac{1}{2}}|u_n|^{\frac{1}{2}})^{**}(2)=0
$$
Thus the inequality \eqref{false} cannot hold in general.
\end{Example}


Let us finally explain why Theorem 1.2 of \cite{FiFoRoSo} required reformulation. For inequalities involving only two norms, like Sobolev or Poincare inequalities, the question about optimality is immediate: we fix one norm and ask how much the second one can be improved. For the Gagliardo--Nirenberg type inequalities the situation is more complicated. Namely, the question about optimality may be asked at least in three ways. 
Considering 
\begin{equation}\label{ME11}
\|\nabla^j u\|_{X}\lesssim\|\nabla^k u\|_{Y}^{\frac{j}{k}}
\|u\|_{Z}^{1-\frac{j}{k}}
\end{equation}
we need to fix two spaces and ask for the optimality of the third one. 

In our main results we follow one possible approach: having two spaces $Y, Z$ corresponding to function and its higher derivative, we constructed/found optimal (roughly speaking) space $X$ corresponding to the middle derivative. 

However, from the point of view of applications, another approach may be desirable. Namely, having spaces $X$ and $Y$ [or $X$ and $Z$] fixed, we need to find the third space $Z$ [or $Y$], possibly optimal, such that (\ref{ME11}) holds.
Such a point of view has already been presented in \cite{FiFoRoSo}. That approach, however, has some disadvantages.

To discuss it, we need to provide new notation. Given a couple of r.i.B.f. spaces $X, Y$ such that $Y\stackrel{\textup{loc}}{\hookrightarrow} X$, we define the space $M(X, Y)$  of pointwise multipliers from $ X $ to $ Y $ using the  following formula
\[
M(X,Y)=\{f\in L^0:fg\in Y{\rm \ for \ each \ }g\in X\}
\]
and equip it with the natural operator norm
$$
\|f\|_{M(X,Y)}:=\sup_{\|g\|_X\leq 1}\|fg\|_Y.
$$
With this norm $M(X,Y)$ becomes ag r.i.B.f. space.
Moreover, the following general version of H\" older inequality follows directly from the definition
$$
\|fg\|_Y\leq\|g\|_{X}\|f\|_{M(X,Y)}.
$$
Notice that $M(X, Y)$ is a kind of generalization of the K\"othe dual of $X$ (i.e. $X'=M(X, L^1)$), but on the other hand, it may be regarded as a kind of point-wise quotient of the space $Y$ by $X$. More information about the spaces of pointwise multipliers may be found in \cite{MP89, KLM13}.

Now we can state an alternative version of Theorem \ref{GBFSR}. Notice that point (i) is actually Theorem 1.2 of \cite{FiFoRoSo}.

\begin{Theorem}\label{GNmult}
Let $j,k\in\mathbb{N}$, $1\leq j<k$. 
\begin{itemize}
\item[(i)]  If $X,Y$ are r.i.B.f. spaces such that 
$$
Y^{\frac{k}{j}}\stackrel{\textup{loc}}{\hookrightarrow} X
$$
and both $Y$ and $M(Y^{k/j},X)^{1-j/k}$ have upper Boyd indices less than $1$, then for all $u\in W^{k,1}_{\loc}(\mathbb{R}^d)$ it holds
\begin{equation}\label{ME21}
\|\nabla^j u\|_{X}\lesssim\|\nabla^k u\|_{Y}^{j/k}
\|u\|_{M(Y^{k/j},X)^{1-j/k}}^{1-j/k}.
\end{equation}

\item[(ii)] If $X,Z$ are r.i.B.f. spaces such that 
$$
Z^{\frac{k}{k-j}}\stackrel{\textup{loc}}{\hookrightarrow} X
$$
and both $Z$ and $M(Z^{k/(k-j)},X)^{j/k}$ have upper Boyd indices less than $1$, then for all $u\in W^{k,1}_{\loc}(\mathbb{R}^d)$ it holds
\begin{equation}\label{ME22}
\|\nabla^j u\|_{X}\lesssim\|\nabla^k u\|_{M(Z^{k/(k-j)},X)^{j/k}}^{j/k}
\|u\|_{Z}^{1-j/k}.
\end{equation}
\end{itemize}
\end{Theorem}
\proof The proof is immediate once we know Theorem \ref{GBFSR}.  
Let us consider point (i). In fact, it is enough to see that 
\[
Y^{j/k}[M(Y^{k/j},X)^{1-j/k}]^{1-j/k}\subset X.
\]
We have by (\ref{propCL1}) 
\[
Y^{j/k}[M(Y^{k/j},X)^{1-j/k}]^{\frac{k-j}{k}}=Y^{k/j}\odot [M(Y^{k/j},X)^{1-j/k}]^{\frac{k}{k-j}}
\]
\[
=Y^{k/j}\odot M(Y^{k/j},X)\subset X. 
\]
\endproof

Notice that the last inclusion above, i.e. $Y^{k/j}\odot M(Y^{k/j}, X)\subset X$, can not be replaced by equality. This is the reason that the formulation of Theorem \ref{GBFSR} is more accurate than the claim of Theorem \ref{GNmult} above (thus also Theorem 1.2 of \cite{FiFoRoSo}). Let us explain it better. We already know that the Gagliardo--Nirenberg inequality 
\[
\|\nabla^j u\|_{X}\lesssim\|\nabla^k u\|_{Y}^{\frac{j}{k}}
\|u\|_{Z}^{1-\frac{j}{k}}.
\]
holds with $X=Y^{j/k}Z^{1-j/k}$ and is optimal in many cases (for more details see  Section \ref{Optimality}). So, to keep potential optimality, having spaces $X, Y$, we need to find the third space $Z$ satisfying $X=Y^{j/k}Z^{1-j/k}$. However, such a space $Z$ may not exist, and the best we can choose is the space 
\[
Z=M(Y^{k/j},X)^{1-j/k}
\]
appearing in the Theorem \ref{GNmult}(i), which for sure gives only inclusion 
\[
Y^{j/k}[M(Y^{k/j},X)^{1-j/k}]^{1-j/k}\subset X,
\] 
but need not give equality. 
Whether this inclusion becomes equality is exactly the problem of (Lozanovskii like) factorization.

Further, we say that $Y$ factorizes $X$ when 
\[
Y\odot M(Y, X)=X.
\]
The origins of factorization come from the Lozanovskii factorization theorem \cite[Theorem 6]{Loz69} (which says that for each B.f. space $X$, there holds $X\odot X'=L^1$), while Pisier has considered factorization in terms of Calderon-Lozanovskii construction in \cite{Pi79} (cf. \cite{Ni85}). For a quite comprehensive study of factorization problems, see, for example, \cite{KLM14, KLM19, Sch10} and references therein. 

In this language, our previous question becomes whether
\[
Y^{k/j}\odot M(Y^{k/j},X)=X? 
\]
This is, whether $Y^{k/j}$ factorizes $X$? If not, then one shouldn't expect optimality of   (\ref{ME21}) (or  (\ref{ME22}), respectively). 

Let us use example to explain the situation. Consider $X=L^{2,\infty}$, $Y=L^{2,1}$ and $j=1,k=2$. Then $Y^2=L^{4,2}$. We have by \cite[Theorem 4]{KLM19}
\[
M(Y^{2},X)=M(L^{4,2},L^{2,\infty})=L^{4,\infty}.
\]
Thus 
\[
Z=M(Y^{2},X)^{1/2}=L^{2,\infty}
\]
and Theorem \ref{GNmult}(i) gives 
\[
\|\nabla u\|_{L^{2,\infty}}\lesssim\|\nabla^2 u\|_{L^{2,1}}^{1/2}
\|u\|_{L^{2,\infty}}^{1/2}.
\]
On the other hand, we know from Theorem \ref{GBFSR} or Corollary \ref{Lorentz} that stronger inequalities hold, i.e. 
\[
\|\nabla u\|_{L^{2,2}}\lesssim\|\nabla^2 u\|_{L^{2,1}}^{1/2}
\|u\|_{L^{2,\infty}}^{1/2}
\]
and 
\[
\|\nabla u\|_{L^{2,\infty}}\lesssim\|\nabla^2 u\|_{L^{2,\infty}}^{1/2}
\|u\|_{L^{2,\infty}}^{1/2}
\]

Concluding, we see that formulation of Theorem \ref{GNmult} may differ from that of Theorem \ref{GBFSR} when respective spaces do not factorize each other. In consequence, only Theorem \ref{GBFSR} has the potential to be optimal and, in fact, it is in many cases, as we will see in the next section.

\section{Optimality}\label{Optimality}

As announced in the introduction, we will show that in the most classical case of $k=2$ and $j=1$, just proved Gagliardo--Nirenberg inequality is optimal for a broad class of r.i.B.f. spaces. 
We have just shown that
\begin{equation}\label{eq:gno}
\|\nabla u\|_{X}\lesssim\|\nabla^2 u\|_{Y}^{\frac{1}{2}}
\|u\|_{Z}^{\frac{1}{2}}
\end{equation}
holds with $X=Y^{1/2}Z^{1/2}$ under some assumptions on $Y,Z$. Now we study the following question of optimality. 
For fixed $Y, Z$, the choice of $X$ is optimal among the r.i.B.f. spaces? More precisely, can $X$ be replaced in \eqref{eq:gno} by strictly smaller r.i.B.f. space $B$ (we mean that $B\subsetneq X$) such that the estimate is still valid for each $u\in  W^{2,1}_{\loc}(\mathbb{R}^d)$? In this section we show that in many cases of spaces $Y,Z$ the choice $X=Y^{1/2}Z^{1/2}$ is optimal.

For the beginning, we need a few simple observations on Calder\'on--Lozanovskii construction.

\begin{Lemma}\label{lemCL2}
Let $Y,Z$ be r.i.B.f. spaces. 
\begin{itemize}
\item[(i)] If $Z\cap L^{\infty}\subset Y\cap L^{\infty}$ then for each $f=f^*\in Y^{1/2}Z^{1/2}$ of the form
\[
f=\sum _{k=1}^{\infty}c_{k}\chi_{[k-1,k)}
\]
there are $g=g^*\in Y$ and $h=h^*\in Z$ of the same form, i.e. 
\[
g=\sum _{k=1}^{\infty}a_{k}\chi_{[k-1,k)},\ \  
h=\sum _{k=1}^{\infty}b_{k}\chi_{[k-1,k)}
\]
such that 
\[
f\leq g^{1/2}h^{1/2}\ \text{ and }\   h\leq g.
\]

\item[(ii)] If $Z\cap L^{1}\subset Y\cap L^{1}$, then for each $f=f^*\in Y^{1/2}Z^{1/2}$ such that $\supp f\subset [0,1]$
there are $g=g^*\in Y$ and $h=h^*\in Z$  satisfying
\[
f\leq g^{1/2}h^{1/2}\ \text{ and }\   h\leq g.
\]
\end{itemize}
\end{Lemma}
\proof (i) 
Let $f=f^*=\sum _{k=1}^{\infty}c_{k}\chi_{[k-1,k)}\in Y^{1/2}Z^{1/2}$. Then there are  $u\in Y,v\in Z$ such that 
\[
f\leq u^{1/2}v^{1/2}.
\]
By \cite[p. 67]{KPS} we have for each $t\geq 0$
\[
f(t)=f^*(t)\leq u^*(t/2)^{1/2}v^*(t/2)^{1/2}.
\]
Denoting $\eta:=D_{1/2}u^*,\gamma:=D_{1/2}v^*$ we get
\[
f\leq \eta^{1/2}\gamma^{1/2},
\]
while $\eta\in Y,\gamma\in Z$ since $Y,Z$ are r.i.B.f. spaces, as the dilation operator is bounded in every r.i.B.f. space.  

Further, we define the averaging operator $T$ (also known as {\it conditional expectation operator})
by the formula
\begin{equation}\label{eq:nondefinedbefore}
T:w\mapsto \sum _{k=1}^{\infty}
\left(\int_{k-1}^kw(t)\,dt \right)\chi _{[k-1,k)}. 
\end{equation}
Since each r.i.B.f. space with the majorant property is the exact interpolation space for the couple $(L^{\infty },L^{1})$, the operator $T$ is non-expanding mapping on such a space (see \cite{BS,KPS}). In consequence, 
\[
r:=T\eta\in Y  \ {\rm\ and\ }\ h:=T\gamma\in Z. 
\]
Both $r$ and $h$ are of the form
\[
r=r^*=\sum _{k=1}^{k}a_{k}\chi_{[k-1,k)},\ \  
h=h^*=\sum _{k=1}^{k}b_{k}\chi_{[k-1,k)}
\]
where 
\[
a_{k}=\int_{k-1}^k\eta(t)dt \ {\rm\ and\ }\ b_k=\int_{k-1}^k\gamma(t)dt. 
\]
We need to see that also $f\leq r^{1/2}h^{1/2}$. Indeed, 
by the H\"{o}lder inequality we get immediately 
$$\begin{aligned}
c_{k}&=\int_{k-1}^kf(t) \,dt\leq \int _{k-1}^k\eta(t)^{1/2}\gamma(t)^{1/2} \,dt \\
&\leq \left(\int_{k-1}^k\eta(t)\,dt\right)^{1/2}\left(\int_{k-1}^k\gamma(t)\,dt\right)^{1/2}=a_k^{1/2}b_k^{1/2}.
\end{aligned}$$
Consequently, it holds $r\in Y\cap L^{\infty}$ and $h\in Z\cap L^{\infty}$, since they are of the special form. Therefore, applying the assumption $Z\cap L^{\infty}\subset Y\cap L^{\infty}$, we see that also 
\[
g:=\max\{r,h\}\in  Y\cap L^{\infty}.
\]
Thus functions $g$ and $h$ have all properties declared in the point (i) since the maximum of nonincreasing functions is also a nonincreasing function. 

(ii) Let $Z\cap L^{1}\subset Y\cap L^{1}$ and choose $f^*=f\in Y^{1/2}Z^{1/2}$ such that $\supp f\subset [0,1]$. Arguing as before, we conclude that there are $\eta=\eta^*\in Y,\gamma=\gamma^*\in Z$ such that 
\[
f\leq \eta^{1/2}\gamma^{1/2}.
\]
Evidently, we can assume that also $\supp\eta\subset[0,1]$ and $\supp\gamma\subset[0,1]$. Then $\eta\in Y\cap L^{1}$ and  $\gamma\in Z\cap L^{1}$. Consequently, thanks to the assumption $Z\cap L^{1}\subset Y\cap L^{1}$, it is enough to take 
\[
g:=\max\{\eta,\gamma\} {\rm \ and\ } h:=\gamma.
\]
\endproof

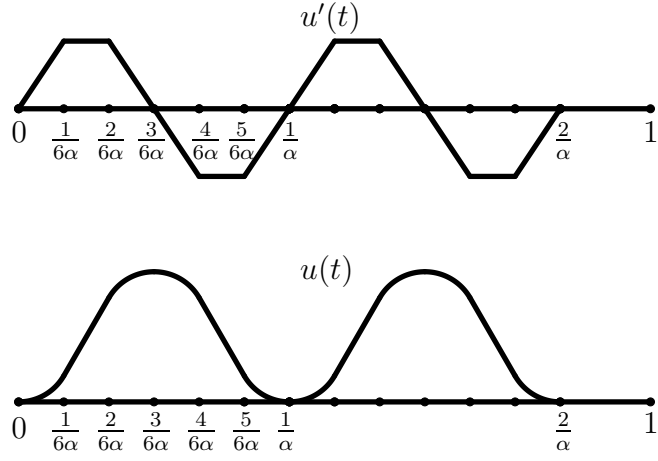
\begin{figure}\label{Figure}
    \centering
    \begin{tikzpicture}[line cap=round,line join=round,>=triangle 45,x=0.6cm,y=0.6cm]
\clip(0,-2) rectangle (17,10);
\draw [line width=2pt] (1,1)-- (15,1);
\draw [line width=2pt] (2,1.5773502691896255)-- (3,3.3094010767585);
\draw [shift={(1,2.1547005383792515)},line width=2pt]  plot[domain=4.71238898038469:5.759586531581287,variable=\t]({1*1.1547005383792515*cos(\t r)+0*1.1547005383792515*sin(\t r)},{0*1.1547005383792515*cos(\t r)+1*1.1547005383792515*sin(\t r)});
\draw [shift={(4,2.732050807568874)},line width=2pt]  plot[domain=4.71238898038469:5.759586531581287,variable=\t]({-1*1.1547005383792515*cos(\t r)+0*1.1547005383792515*sin(\t r)},{0*1.1547005383792515*cos(\t r)+-1*1.1547005383792515*sin(\t r)});
\draw [shift={(4,2.732050807568874)},line width=2pt]  plot[domain=4.71238898038469:5.759586531581287,variable=\t]({1*1.1547005383792515*cos(\t r)+0*1.1547005383792515*sin(\t r)},{0*1.1547005383792515*cos(\t r)+-1*1.1547005383792515*sin(\t r)});
\draw [shift={(7,2.154700538379251)},line width=2pt]  plot[domain=4.71238898038469:5.759586531581287,variable=\t]({-1*1.1547005383792515*cos(\t r)+0*1.1547005383792515*sin(\t r)},{0*1.1547005383792515*cos(\t r)+1*1.1547005383792515*sin(\t r)});
\draw [line width=2pt] (6,1.5773502691896253)-- (5,3.3094010767585);
\draw [shift={(13,2.1547005383792506)},line width=2pt]  plot[domain=4.71238898038469:5.759586531581287,variable=\t]({-1*1.1547005383792515*cos(\t r)+0*1.1547005383792515*sin(\t r)},{0*1.1547005383792515*cos(\t r)+1*1.1547005383792515*sin(\t r)});
\draw [shift={(10,2.7320508075688736)},line width=2pt]  plot[domain=4.71238898038469:5.759586531581287,variable=\t]({1*1.1547005383792515*cos(\t r)+0*1.1547005383792515*sin(\t r)},{0*1.1547005383792515*cos(\t r)+-1*1.1547005383792515*sin(\t r)});
\draw [shift={(10,2.7320508075688736)},line width=2pt]  plot[domain=4.71238898038469:5.759586531581287,variable=\t]({-1*1.1547005383792515*cos(\t r)+0*1.1547005383792515*sin(\t r)},{0*1.1547005383792515*cos(\t r)+-1*1.1547005383792515*sin(\t r)});
\draw [shift={(7,2.154700538379251)},line width=2pt]  plot[domain=4.71238898038469:5.759586531581287,variable=\t]({1*1.1547005383792515*cos(\t r)+0*1.1547005383792515*sin(\t r)},{0*1.1547005383792515*cos(\t r)+1*1.1547005383792515*sin(\t r)});
\draw [line width=2pt] (12,1.5773502691896248)-- (11,3.3094010767584994);
\draw [line width=2pt] (8,1.5773502691896253)-- (9,3.3094010767585);
\draw [line width=2pt] (1,7.5)-- (15,7.5);
\draw [line width=2pt] (1,7.5)-- (2,9);
\draw [line width=2pt] (2,9)-- (3,9);
\draw [line width=2pt] (3,9)-- (4,7.5);
\draw [line width=2pt] (7,7.5)-- (6,6);
\draw [line width=2pt] (6,6)-- (5,6);
\draw [line width=2pt] (5,6)-- (4,7.5);
\draw [line width=2pt] (13,7.5)-- (12,6);
\draw [line width=2pt] (12,6)-- (11,6);
\draw [line width=2pt] (11,6)-- (10,7.5);
\draw [line width=2pt] (7,7.5)-- (8,9);
\draw [line width=2pt] (8,9)-- (9,9);
\draw [line width=2pt] (9,9)-- (10,7.5);
\draw (7,10.1) node[anchor=north west] {$u'(t)$};
\draw (7,4.5) node[anchor=north west] {$u(t)$};
\draw (0.6,7.45) node[anchor=north west] {$0$};
\draw (0.6,0.9) node[anchor=north west] {$0$};
\draw (1.45,7.45) node[anchor=north west] {$\frac{1}{6\alpha}$};
\draw (2.45,7.45) node[anchor=north west] {$\frac{2}{6\alpha}$};
\draw (3.35,7.45) node[anchor=north west] {$\frac{3}{6\alpha}$};
\draw (5.35,7.45) node[anchor=north west] {$\frac{5}{6\alpha}$};
\draw (4.55,7.45) node[anchor=north west] {$\frac{4}{6\alpha}$};
\draw (6.55,7.45) node[anchor=north west] {$\frac{1}{\alpha}$};
\draw (12.6,7.45) node[anchor=north west] {$\frac{2}{\alpha}$};
\draw (14.6,7.45) node[anchor=north west] {$1$};
\draw (1.45,1) node[anchor=north west] {$\frac{1}{6\alpha}$};
\draw (2.45,1) node[anchor=north west] {$\frac{2}{6\alpha}$};
\draw (3.45,1) node[anchor=north west] {$\frac{3}{6\alpha}$};
\draw (4.45,1) node[anchor=north west] {$\frac{4}{6\alpha}$};
\draw (5.45,1) node[anchor=north west] {$\frac{5}{6\alpha}$};
\draw (6.45,1) node[anchor=north west] {$\frac{1}{\alpha}$};
\draw (12.6,1) node[anchor=north west] {$\frac{2}{\alpha}$};
\draw (14.6,1) node[anchor=north west] {$1$};
\begin{scriptsize}
\draw [fill=black] (1,1) circle (1.5pt);
\draw [fill=black] (7,1) circle (1.5pt);
\draw [fill=black] (13,1) circle (1.5pt);
\draw [fill=black] (15,1) circle (1.5pt);
\draw [fill=black] (1,1) circle (1.5pt);
\draw [fill=black] (3,1) circle (1.5pt);
\draw [fill=black] (2,1) circle (1.5pt);
\draw [fill=black] (4,1) circle (1.5pt);
\draw [fill=black] (5,1) circle (1.5pt);
\draw [fill=black] (6,1) circle (1.5pt);
\draw [fill=black] (7,1) circle (1.5pt);
\draw [fill=black] (6,1) circle (1.5pt);
\draw [fill=black] (7,1) circle (1.5pt);
\draw [fill=black] (13,1) circle (1.5pt);
\draw [fill=black] (11,1) circle (1.5pt);
\draw [fill=black] (12,1) circle (1.5pt);
\draw [fill=black] (10,1) circle (1.5pt);
\draw [fill=black] (9,1) circle (1.5pt);
\draw [fill=black] (8,1) circle (1.5pt);
\draw [fill=black] (7,1) circle (1.5pt);
\draw [fill=black] (8,1) circle (1.5pt);
\draw [fill=black] (1,7.5) circle (1.5pt);
\draw [fill=black] (1,7.5) circle (1.5pt);
\draw [fill=black] (15,7.5) circle (1.5pt);
\draw [fill=black] (7,7.5) circle (1.5pt);
\draw [fill=black] (13,7.5) circle (1.5pt);
\draw [fill=black] (1,7.5) circle (1.5pt);
\draw [fill=black] (3,7.5) circle (1.5pt);
\draw [fill=black] (2,7.5) circle (1.5pt);
\draw [fill=black] (4,7.5) circle (1.5pt);
\draw [fill=black] (5,7.5) circle (1.5pt);
\draw [fill=black] (6,7.5) circle (1.5pt);
\draw [fill=black] (7,7.5) circle (1.5pt);
\draw [fill=black] (6,7.5) circle (1.5pt);
\draw [fill=black] (7,7.5) circle (1.5pt);
\draw [fill=black] (13,7.5) circle (1.5pt);
\draw [fill=black] (11,7.5) circle (1.5pt);
\draw [fill=black] (12,7.5) circle (1.5pt);
\draw [fill=black] (10,7.5) circle (1.5pt);
\draw [fill=black] (9,7.5) circle (1.5pt);
\draw [fill=black] (8,7.5) circle (1.5pt);
\draw [fill=black] (7,7.5) circle (1.5pt);
\draw [fill=black] (8,7.5) circle (1.5pt);
\draw [fill=black] (7,7.5) circle (1.5pt);
\draw [fill=black] (4,7.5) circle (1.5pt);
\draw [fill=black] (7,7.5) circle (1.5pt);
\draw [fill=black] (7,7.5) circle (1.5pt);
\draw [fill=black] (5,7.5) circle (1.5pt);
\draw [fill=black] (6,7.5) circle (1.5pt);
\draw [fill=black] (13,7.5) circle (1.5pt);
\draw [fill=black] (10,7.5) circle (1.5pt);
\draw [fill=black] (7,7.5) circle (1.5pt);
\draw [fill=black] (10,7.5) circle (1.5pt);
\draw [fill=black] (13,7.5) circle (1.5pt);
\draw [fill=black] (13,7.5) circle (1.5pt);
\draw [fill=black] (11,7.5) circle (1.5pt);
\draw [fill=black] (12,7.5) circle (1.5pt);
\draw [fill=black] (9,7.5) circle (1.5pt);
\draw [fill=black] (8,7.5) circle (1.5pt);
\draw [fill=black] (7,7.5) circle (1.5pt);
\draw [fill=black] (8,7.5) circle (1.5pt);
\draw [fill=black] (7,7.5) circle (1.5pt);
\draw [fill=black] (7,7.5) circle (1.5pt);
\draw [fill=black] (7,7.5) circle (1.5pt);
\draw [fill=black] (7,7.5) circle (1.5pt);
\draw [fill=black] (9,7.5) circle (1.5pt);
\draw [fill=black] (8,7.5) circle (1.5pt);
\end{scriptsize}
\end{tikzpicture}
    \caption{Sketch of the construction of $u'$ and $u$ in Lemma \ref{lan}}
    \label{fig:my_label}
\end{figure}

\begin{Lemma}\label{lan}
For every $\alpha\geq1$ there exists $u\in W^{2,\infty}(\mathbb{R})$, $\supp u=[0,1]$, such that
\begin{enumerate}[\upshape(i)]
\item 
$
\bigl|\left\{|u'|\geq \alpha\right\}\bigr|\geq 1/6
$
\item $|u(t)|\leq 1/3$ and $|u''(t)|\leq 6\alpha^2$.
\end{enumerate}
\end{Lemma}
\proof[Proof of Lemma]
On $(0,1/\alpha)$ we define a continuous piecewise affine $u'(t)$  as
$$
u'(t):=
 \begin{cases}
6\alpha^2t  &\mbox{if  } t\in[0,\frac{1}{6\alpha}]\\
 \alpha &\mbox{if } t\in[\frac{1}{6\alpha},\frac{2}{6\alpha}] \\
3\alpha-6\alpha^2t&\mbox{if } t\in[\frac{2}{6\alpha},\frac{3}{4\alpha}] \\
-\alpha &\mbox{if } t\in[\frac{4}{6\alpha},\frac{4}{5\alpha}]  \\
-5\alpha+6\alpha^2t  &\mbox{if } t\in[\frac{5}{6\alpha},\frac{1}{\alpha}]. \\
  \end{cases}
$$ 
and $u(t)=\int_{0}^tu'(s)\,ds$.
This formula can be copied several times on intervals $[k/\alpha,(k+1)/\alpha]$ for $k\in\mathbb{N}, k<\alpha$ by $u(t):=u(t-1/\alpha)$. If $\alpha$ is integer number, this will cover whole $[0,1]$ interval, otherwise we extend by constant zero, as we illustrate by Figure \ref{Figure} for case $\alpha\in(2,3)$. Note that for $\alpha$ integer, we may get constant 1/3 in (i), but for $\alpha$ non-integer the achieved constant has to be lowered to 1/6.
\endproof


We are finally ready to prove the main theorem of the paper. 
\proof[Proof of Theorem \ref{optimal}]
(i) Assume that  $Z\cap L^{\infty}\subset Y\cap L^{\infty}$ and suppose $X\cap L^{\infty}\not \subset B\cap L^{\infty}$. We will construct a function $\eta\in  W^{2,1}_{\loc}(\mathbb{R}^d)$ such that 
\[
\|\nabla^2\eta\|_Y^{\frac{1}{2}}\|\eta\|_Z^{\frac{1}{2}}< \infty,
\]
but
\[
\|\nabla \eta\|_B=\infty.  
\]

Firstly we will consider the case of $d=1$. 
Since $X\cap L^{\infty}\not \subset B\cap L^{\infty}$ there is $\tilde{f}=(\tilde{f})^*\in  (X\cap L^{\infty})\backslash (B\cap L^{\infty})$ with $\|\tilde{f}\|_X=1$ and observe that $\tilde{f}$ can be chosen in the form of
\[
\tilde{f}=\sum_{k=1}^{\infty}\tilde{c_k}\chi_{[k-1,k)}.
\]

On the other hand, $X=Y^{1/2}Z^{1/2}$, thus from Lemma \ref{lemCL2} it follows that there are $g\in Y$, $h\in Z$ such that $\tilde{f}\leq g^{1/2}h^{1/2}$,
\begin{equation}\label{eq:optchoice}
g=g^*=\sum_{k=1}^{\infty}a_k\chi_{[k-1,k)},\ \ 
\\
h=h^*=\sum_{k=1}^{\infty}b_k\chi_{[k-1,k)}
\end{equation}
and $h\leq g$. Define 
$$
f=g^{1/2}h^{1/2}.
$$
Then $f=f^*$ is of the form
\[
f=\sum_{k=1}^{\infty}c_k\chi_{[k-1,k)}
\]
and also $f\in X\backslash B$, since $0\leq \tilde{f}\leq f$. 
In particular,
$\|f\|_B=\infty$.

By definition of $f$ we have $c_k=\sqrt{a_kb_k}$ for each $k$, thus  for each $k$
\[
b_k\leq c_k\leq a_k,
\]
since $h\leq g$.
Now for each $k=1,2,...$ we define $u_k$ as $u$ from Lemma \ref{lan} applied with $\alpha=c_k/b_k$. Further we define 
\[
\eta(t)=\sum_{k=1}^{\infty}b_ku_k(t-k+1).
\]
In consequence,
\[
\eta'(t)=\sum_{k=1}^{\infty}b_ku'_k(t-k+1)
\]
From Lemma \ref{lan} it follows
$$
\left|\left\{|b_k u'_k(t)|\geq c_k\right\}\right|\geq \frac{1}{6},
$$
which implies that for each $n\in\mathbb{N}$
\begin{equation}\label{eq:opt1}
(\eta')^{*}\geq \sum_{k=1}^{n} c_k\chi_{[\frac{k-1}{6},\frac{k}{6})}=D_{6}\left(\sum_{k=1}^{n} c_k\chi_{[k-1,k)}\right).
\end{equation}
Since dilation operators are bounded on r.i.B.f. spaces we get 
\begin{equation}\label{eq:opt2}\begin{aligned}
\left\|D_{6}\left(\sum_{k=1}^{n} c_k\chi_{[k-1,k)}\right)\right\|_B&\geq 
\|D_{1/6}\|_{B\to B}^{-1}\left\|D_{1/6}D_6\left(\sum_{k=1}^{n} c_k\chi_{[k-1,k)}\right)\right\|_B
\\
&=\|D_{1/6}\|_{B\to B}^{-1}\left\|\sum_{k=1}^{n} c_k\chi_{[k-1,k)}\right\|_B.
\end{aligned}
\end{equation}
The limit for $n\to\infty$ of the term on right-hand side is $\infty$ since $\sum_{k=1}^{n} c_k\chi_{[k-1,k)}\to f$ pointwise, while $B$ has the Fatou property and $f\not \in B$. Together with estimates \eqref{eq:opt1} and \eqref{eq:opt2} it implies that 
\[
\|\eta'\|_B=\infty.
\]
Furthermore, 
\[
\eta''(t)=\sum_{k=1}^{\infty}b_ku''_k(t-k+1)
\]
and the way we have chosen functions $u_k$ ensures that for each $k=1,2,...$ and each $t\in[k-1,k)$
$$
b_k|u''_k(t-k+1)|\leq 6b_k\left(\frac{c_k}{b_k}\right)^2\leq 6 a_k,
$$ 
while 
\[
u''_k(t-k+1)=0 {\rm \ for\ each \ }t\notin [k-1,k).
\]
Using this and \eqref{eq:optchoice} we get
\[
|\eta''|\leq \sum_{k=1}^{\infty}a_k\chi_{[k-1,k)}=g.
\]
In consequence $\|\eta''\|_Y<\|g\|_Y<\infty$. By \eqref{eq:optchoice} we have $\|\eta\|_Z\leq \|h\|_Z<\infty$, since 
\[
\eta\leq \sum_{k=1}^{\infty}b_k\chi_{[k-1,k)}=h.
\]
This finishes the proof of point (i) in the case of $d=1$. 

In case $d>1$ we proceed in the same fashion. Firstly we select $f,g,h$ together with $b_k,c_k,a_k$ and define functions $u_k$ exactly as above. Further, define $w'$ to be a piecewise affine continuous function with slots in points $(0,0),(1/4,1),(3/4,-1),(1,0)$. Precisely the formula is
$$
w'(t):=
 \begin{cases}
4t  &\mbox{if  } t\in[0,\frac{1}{4}]\\
2-4t &\mbox{if } t\in[\frac{1}{4},\frac{3}{4}] \\
4t-4&\mbox{if } t\in[\frac{3}{4},1] \\
0 &\mbox{else. }
  \end{cases}
$$ 
We set for $t>0$
\[
w(t)=\int_0^tw'(s)\,ds.
\]
Then $w\in W^{2,\infty}(\mathbb{R})$ has the following properties:\\
(a) $\supp w= [0,1]$\\
(b) $|w(t)|\leq 1/4$, $|w'(t)|\leq 1$ and $|w''(t)|\leq 4$ for each $t\in [0,1]$\\
(c) $|w(t)|\geq 1/8$ for each $t\in [1/4,3/4]$.

Finally, we are ready to define the desired function $\eta$ on $\mathbb{R}^d$. We put
\[
\eta=\sum_{k=1}^{\infty}b_kv_k,
\]
where for each $k=1,2,3,...$
\[
v_k(x_1,...,x_d)=u_k(x_1-k+1)w(x_2)w(x_3)...w(x_{d}).
\]
It remains to estimate the function $\eta$ and its derivatives analogously as in the previous part of the proof. By Lemma \ref{lan} and points (a), (b) above we have 
\[
\eta^*\leq h.
\] 
Moreover, notice that $|u_k'|\leq \frac{c_k}{b_k}$, therefore for each fixed $1<i\leq d$
\[
|b_ku_k'(x_1-k+1)w'(x_i)\Pi_{j=2,j\not=i}^{d}w(x_j)|\leq b_k\frac{c_k}{b_k}\leq a_k,
\]
while respective estimates for another second order derivatives appearing in $\nabla^2\eta$ are immediate.
In consequence,
\[
(\nabla^2\eta)^*\leq 64^dd^2\sum_{k=1}^{\infty}a_k\chi_{[k-1,k)}=64^dd^2g.
\]
Finally, point (c) above implies that 
\[
|\nabla \eta|\geq \left|\frac{\partial \eta}{\partial x_1}\right| \geq \sum_{k=1}^{n}\frac{1}{8^d}c_k\chi_{C_k},
\]
 where $C_k=\left\{|b_k u'_k|\geq c_k\right\}\times [1/4,3/4]^{d-1}$.
Thus $|C_k|\geq \frac{1}{2^{d+2}}$ by definition of $u_k$'s. 
In consequence, for each
$n\in\mathbb{N}$
\[
|\nabla \eta|^{*}\geq \frac{1}{8^d}\sum_{k=1}^{n} c_k\chi_{[\frac{k-1}{2^{d+2}},\frac{k}{2^{d+2}})}=\frac{1}{8^d}D_{2^{d+2}}\left(\sum_{k=1}^{n} c_k\chi_{[k-1,k)}\right)
\]
and we can finish the proof as before. 
\vspace{6mm}

(ii) Consider the case $Z\cap L^{1} \subset Y\cap L^{1}$  and suppose  $X\cap L^{1}\not \subset B\cap L^{1}$. 
This time we will construct a sequence $(\eta_n)\subset W^{2,1}_{\loc}(\mathbb{R}^d)$ such that 
\[
\frac{\|\nabla \eta_n\|_{B}}{\|\nabla^2 \eta_n\|_{Y}^{\frac{1}{2}}
\|\eta_n\|_{Z}^{\frac{1}{2}}}\to \infty.
\] 
 Once again we start with the case of $d=1$. It follows that there is
 $\tilde{f}=(\tilde{f})^*\in X\backslash B$ such that
\[
\supp {\tilde{f}}=[0,1].
\]
On the other hand, $\tilde{f}\in X=Y^{1/2}Z^{1/2}$ so by Lemma \ref{lemCL2} there are $g^*=g\in Y$, $h^*=h\in Z$ such that $\tilde{f}\leq g^{1/2}h^{1/2}$ and $h\leq g$. Define for each $n$
\[
g_n=g_n^*=\sum_{k=1}^{\infty}g\left(\frac{k}{n}\right)\chi_{[\frac{k-1}{n},\frac{k}{n}]},
\]
\[
h_n=h_n^*=\sum_{k=1}^{\infty}h\left(\frac{k}{n}\right)\chi_{[\frac{k-1}{n},\frac{k}{n}]}.
\]
Set
$$
a_{n,k}:=g\left(\frac{k}{n}\right), b_{n,k}:=h\left(\frac{k}{n}\right)\textup{ and } c_{n,k}:=\sqrt{a_{n,k}b_{n,k}}.
$$
Note that since $0\leq h_n\leq h_{n+1}\leq h$ and $0\leq g_n\leq g_{n+1}\leq g$ by the lattice property we obtain that 
\[
\|h_n\|_{Z}\leq \|h\|_{Z} {\rm \ and\ } \|g_n\|_{Y}\leq \|g\|_{Y}. 
\]

Further, we put $f_n=g_n^{1/2}h_n^{1/2}$. 
Then $f_n$ is of the form
\[
f_n=\sum_{k=1}^{\infty}c_{n,k}\chi_{[\frac{k-1}{n},\frac{k}{n}]}.
\]
The pointwise monotone convergence $0\leq f_n\uparrow f= g^{1/2}h^{1/2}$ combined with the Fatou property of $B$ yields 
 \[
 \left\|f_n\right\|_B\to\|f\|_B\geq \|\tilde{f}\|_B=\infty.
 \]
Moreover, $c_{n,k}=\sqrt{a_{n,k}b_{n,k}}$  and  
\[
b_{n,k}\leq c_{n,k}\leq a_{n,k} {\rm \ for\ each \ }k=1,2,...,n.
\]
Now for each $n$ and $k=1,2,...,n$ we select $u_{n,k}$ from Lemma  \ref{lan} applied with $c=c_{n,k}/b_{n,k}$. Finally we define 
\[
\eta_n(t)=\sum_{k=1}^{n}b_{n,k}u_{n,k}(nt-k+1).
\]
Then
\[
\eta'_n(t)=n\sum_{k=1}^{n}b_{n,k}(u_{n,k})'(nt-k+1)
\]
and we have by Lemma \ref{lan}
$$
\left|\left\{b_{n,k} |u_{n,k}'|\geq c_{n,k}\right\}\right|\geq \frac{1}{6n}.
$$
Hence
\[
(\eta'_n)^{*}(t)\geq n \sum_{k=1}^{n} c_{n,k}\chi_{[\frac{k-1}{6n},\frac{k}{6n})}=nD_{6}\left(\sum_{k=1}^{n} c_{n,k}\chi_{[\frac{k-1}{n},\frac{k}{n}]}\right).
\]
We have
$$
n\left\|\sum_{k=1}^{n}  c_{n,k}\chi_{[\frac{k-1}{n},\frac{k}{n})}\right\|_B=n \left\|D_{1/6}D_6\left(\sum_{k=1}^{n}  c_{n,k}\chi_{[\frac{k-1}{n},\frac{k}{n})}\right)\right\|_B\leq \|D_{1/6}\|_{B\to B}\|\eta'_n(t)\|_B,
$$
since dilations are bounded on $B$. 

On the other hand, we have 
\[
|\eta_n(t)|\leq \sum_{k=1}^{n}b_{n,k}\chi_{[\frac{k-1}{n},\frac{k}{n})}\leq h_n
\]
and so $\|\eta_n\|_Z\leq\|h_n\|_Z\leq \|h\|_Z<\infty$. Moreover, 
\[
\eta''_n(t)=n^2\sum_{k=1}^{n}b_{n,k}u_{n,k}''(nt-k+1)
\]
and, since 
$$
b_{n,k}|u''_{n,k}(nt-k+1)|\leq b_{n,k}\left(\frac{c_{n,k}}{b_{n,k}}\right)^2=a_{n,k},
$$ 
we get
\[
|\eta''_n(t)|\leq n^2\sum_{k=1}^{n}a_{n,k}\chi_{[\frac{k-1}{n},\frac{k}{n})}=n^2g_n
\]
and consequently $\|\eta''_n\|_Y\leq n^2\|g_n\|_Y\leq n^2\|g\|_Y$. Finally, we have 
\[
\frac{\|\eta'_n\|_{B}}{\|\eta''_n\|_{Y}^{\frac{1}{2}}
\|\eta_n\|_{Z}^{\frac{1}{2}}}\geq 
\frac{n\|f_n\|_{B}}{(n^2\|g_n\|_{Y})^{\frac{1}{2}}
\|h_n\|_{Z}^{\frac{1}{2}}}\to \infty
\]
and the proof of point (ii) is finished in the case of $d=1$. When $d>1$ the proof is a mixture of the above argument together with the idea used for the case $d>1$ in the proof of point (i). More precisely, we define 
\[
\eta_n=\sum_{k=1}^{n}b_{n,k}v_{n,k},
\]
where 
\[
v_{n,k}(x_1,x_2,...,x_d)=u_{n,k}(nx_1-k+1)w(x_2)...w(x_d)
\]
and the rest of the proof looks as before.
\endproof

Considering the spaces on the right-hand side of the Gagliardo--Nirenberg inequality satisfy assumption $Z\subset Y$, then Theorem \ref{optimal} takes the following simplified form.

\begin{Corollary}\label{full}
Let $Y,Z$ be r.i.B.f. spaces satisfying $Z\subset Y$. If 
the Gagliardo--Nirenberg inequality
\begin{equation*}
\|\nabla u\|_{B}\lesssim\|\nabla^2 u\|_Y^{\frac{1}{2}}\|u\|_Z^{\frac{1}{2}}
\end{equation*}
holds for all $u\in W^{2,0}(\mathbb{R}^d)$, then $Y^{1/2}Z^{1/2}\subset B$.
\end{Corollary}
\proof
We need only to see that if $Z\subset Y$ and $Y^{1/2}Z^{1/2}\not \subset B$ then either assumptions of point (i) or point (ii) of Theorem \ref{optimal} are satisfied. First of all, notice that $Z\subset Y$ implies both 
 \[
 Z\cap L^{\infty}\subset Y\cap L^{\infty}{\ \text{and}\ }  Z\cap L^{1}\subset Y\cap L^{1}.
 \]
Now let $0\leq f\in Y^{1/2}Z^{1/2}\backslash B$ and define 
\[
f_a=f\chi_{\{f\geq 1\}}{\ \text{and}\ }  f_b=f\chi_{\{f< 1\}}.
\]
Then either $f_b$ belongs to $(Y^{1/2}Z^{1/2})\backslash B$  or $f_a$ belongs to $(Y^{1/2}Z^{1/2})\backslash B$ (or both are valid). In case  $f_b\in (Y^{1/2}Z^{1/2})\backslash B$ we apply Theorem \ref{optimal} (i), since $f_b\in L^{\infty}$. Otherwise, we apply Theorem \ref{optimal} (ii), since $f_a\in L^{1}$.
\endproof

Let us apply previous thoughts on the Orlicz spaces. As mentioned in the introduction, Ka{\l}amajska and Pietruska-Pa{\l}uba \cite{KP} studied the  Gagliardo--Nirenberg inequality 
\[
\|\nabla u\|_{\varphi}\lesssim\|\nabla^2 u\|_{\varphi_1}^{\frac{1}{2}}\|u\|_{\varphi_2}^{\frac{1}{2}}.
\]
They proved the choice of $L^{\varphi}$ for $\varphi$ satisfying
\[
\varphi^{-1}\approx \sqrt{\varphi_1^{-1}\varphi_2^{-1}}
\]
is valid when spaces satisfy some additional technical assumptions. Later the result was relaxed of some of these assumptions and the optimality of the choice $X=L^{\varphi}$ among all Orlicz spaces was given in \cite{FiFoRoSo}. Theorem \ref{optimal} explains that it is already the optimal choice of space among all r.i.B.f. spaces, provided $L^{\varphi_2}\subset L^{\varphi_1}$.

Concluding, we see that the choice of $X=Y^{1/2}Z^{1/2}$ is optimal among all r.i.B.f. spaces in the Gagliardo--Nirenberg inequality \eqref{ME}, provided that $Z\subset Y$. However, the assumption $Z\subset Y$ is quite restrictive and does not apply to the most classical r.i.B.f. spaces on  $\mathbb{R}_+$ (Lebesgue spaces, Lorentz spaces, etc.), since usually there is no inclusion between such spaces on  $\mathbb{R}_+$. It appears, however, that manoeuvring between points (i) and  (ii) of Theorem \ref{optimal} we can use it to give an almost complete answer to the question about the optimality of (\ref{ME}) among Lorentz spaces posted in \cite{FiFoRoSo} formulated as Corollary \ref{Lorentz}. 

\proof[Proof of Corollary \ref{Lorentz}]
The choice of $P=(1/(2Q)+1/(2R))^{-1}$ is the only possible. This is a consequence of the scaling argument (see \cite[Theorem 1.1]{FiFoRoSo} for details) and of the shape of fundamental functions of Lorentz spaces. Thus we need only to explain that 
\[
\|\nabla u\|_{P,p'}\lesssim\|\nabla^2 u\|_{Q,q}^{\frac{1}{2}}\|u\|_{R,r}^{\frac{1}{2}}
\]
does not hold when $p'<p:=(1/(2q)+1/(2r))^{-1}$ for $P=(1/(2Q)+1/(2R))^{-1}$.
We will consider three cases. Firstly, if $R<Q$ then 
$$
L^{R,r}\cap L^{\infty}\subset  L^{Q,q}\cap L^{\infty}. 
$$
Thus we can  apply point (i) of Theorem \ref{optimal} with $B=L^{P,p'}$, since $L^{P,p'}\subsetneq L^{P,p}$.

In the second case, when $R>Q$, we have
$$
L^{R,r}\cap L^{1}\subset L^{Q,q}\cap L^{1}
$$
and we apply point (ii) of Theorem \ref{optimal}.

Finally, when $R=Q$ and $r<q$, then assumptions of Corollary \ref{full} are satisfied.
\endproof


\section*{Acknowledgments}

The first named author was supported by the National Science Center (Narodowe Centrum Nauki), Poland (project
no.~2017/26/D/ST1/00060). 

The second and the third named authors were supported by the grant number GJ20-19018Y of the Grant Agency of the Czech Republic.

We would like to thank our colleagues and friends for fruitful discussions on the topic, especially to Ale\v{s} Nekvinda and our late friend and mentor Jan Mal\'{y}, whose memory we wish to dedicate this article.

\bibliographystyle{plain}
\bibliography{clanky}
\end{document}